\documentclass{article}
\usepackage{latexsym}
\usepackage{amssymb}
\usepackage{fullpage,amsthm, amsfonts, amsmath, stmaryrd}

\def\0{{\bf 0}}
\def\1{{\bf 1}}

\def\j{{\hat{\jmath}}}
\def\t{{^{\top\!\!\!}}}
\def\stp{\bowtie}

\theoremstyle{plain}
\newtheorem{thm}{Theorem}
\newtheorem{lem}[thm]{Lemma}
\newtheorem{prp}[thm]{Proposition}

\theoremstyle{definition}

\newcommand{\Zz}{{\mathbb Z}}
\newcommand{\adj}{{\,\sim\,}}

\begin{document}
\title{Godsil-McKay Switching and Isomorphism}
\author{Aida Abiad$^*$, Andries E. Brouwer, Willem H. Haemers$^*$
\\
{\small $^*$Tilburg University, Tilburg, The Netherlands}
\\
{\small (e-mails: {\tt A.AbiadMonge@uvt.nl}, {\tt aeb@cwi.nl}, {\tt haemers@uvt.nl})
}
}
\date{}
\maketitle
\abstract{
\noindent
Godsil-McKay switching is an operation on graphs that doesn't change the spectrum of the adjacency matrix.
Usually (but not always) the obtained graph is non-isomorphic with the original graph.
We present a straightforward sufficient condition for being isomorphic after switching,
and give examples which show that this condition is not necessary.
For some graph products we obtain sufficient conditions for being non-isomorphic after switching.
As an example we find that the tensor product of the $\ell\times m$ grid ($\ell>m\geq 2$) and
a graph with at least one vertex of degree two is not determined by its adjacency spectrum.
}
\\[3pt]

\section{Introduction}

An important activity in algebraic graph theory is to decide
if a graph is determined by the spectrum of the adjacency matrix
(see the surveys \cite{DH2003,DH2009}).
Godsil-McKay switching is an operation on a graph that does not change the
spectrum, and provides a tool for disproving existence of such a characterization.
For this operation to work the graph needs a special structure.
However, the presence of this structure doesn't imply that
the graph is not determined by its spectrum; it may be that after switching the graph is isomorphic with the original one.
In this note we investigate this phenomenon.
We hoped to find some useful criteria for isomorphism after switching.
Unfortunately we found some strange examples, which indicate that there is not much hope for
such a criterium.
Instead we obtain some necessary conditions and show how they can be used
to guarantee non-isomorphism after switching for some graph products.

\section{Godsil-McKay switching}\label{GM}

Two graphs with the same (adjacency) spectrum are called {\em cospectral}.
Godsil and McKay \cite{GM1982} introduced the following construction method for cospectral graphs.
\begin{prp}\label{GMswitching}
Let $G$ be a graph and let $\{X_1,\ldots,X_\ell,Y\}$ be a
partition of the vertex set $V(G)$ of $G$.
Suppose that for every vertex $x\in Y$ and every $i\in\{1,\ldots,\ell\}$, $x$ has either $0$,
$\frac{1}{2}|X_i|$ or $|X_i|$ neighbors in $X_i$.
Moreover, suppose that for all $i,j\in\{1,\ldots,\ell\}$ every vertex $x\in X_i$ has the same number of neighbors in $X_j$.
Make a new graph $G'$ as follows.
For each $x\in Y$ and $i\in\{1,\ldots,\ell\}$ such that $x$ has $\frac{1}{2}|X_i|$ neighbors
in $X_i$ delete the corresponding $\frac{1}{2}|X_i|$ edges and join $x$ instead to the
$\frac{1}{2}|X_i|$ other vertices in $X_i$.
Then $G$ and $G'$ are cospectral.
\end{prp}
\noindent
The operation that changes $G$ into $G'$ is called \emph{Godsil-McKay switching},
and the considered partition is a \emph{(Godsil-McKay) switching partition}.
In many applications $\ell=1$.
Then the above condition requires that $X=X_1$ induces a regular subgraph of $G$,
and that each vertex in $Y$ has $0$, $\frac{1}{2}|X|$ or $|X|$ neighbors in $X$.
Such a set $X$ will be called a \emph{(Godsil-McKay) switching set}.
In this note we look for conditions
on a switching set under which $G$ and $G'$ are isomorphic.

Let $G$ be a graph with adjacency matrix $A$ and switching set $X$.
Let $B$ be the submatrix of $A$ corresponding to $X$.
Then
\[
A=\left[
\begin{array}{cc}
B & M \\
M^\t & C
\end{array}
\right],
\ \ \mbox{with}
\ \ M=\left[
\begin{array}{ccc}
N & J & O
\end{array}
\right],
\]
where $BJ=kJ$ for some $k\in\{0,\ldots,|X|-1\}$, and $N^\t J = \frac{1}{2}|X| J$.
Note that not every (but at least one) type of block $N$, $J$ or $O$ needs to be present.
Let $G'$ be the graph with adjacency matrix $A'$ obtained by Godsil-McKay switching with respect
to $X$ in $G$.
Then
\[
A'=\left[
\begin{array}{cc}
B & M' \\
M'^\t & C
\end{array}
\right],
\ \ \mbox{with}
\ \ M'=\left[
\begin{array}{ccc}
J-N & J & O
\end{array}
\right].
\]
With the above notation, the following proposition is straightforward.
\begin{prp}\label{isoGM}
If there exist permutation matrices $P$ and $Q$ such
that $PBP^{\t}=B$, $PMQ^{\t}=M'$ and $QCQ^{\t}=C$, then
$G$ and $G'$ are isomorphic.
\end{prp}
\noindent
Any pair of vertices in $G$ is a switching set, but such a set always
satisfies the above proposition, so switching produces isomorphic graphs.
However, if $|X|\geq 4$ then Proposition~\ref{isoGM} is in not automatically
satisfied and Godsil-McKay switching usually (but not always) produces
non-isomorphic graphs. To prove that $G$ and $G'$ are non-isomorphic
it would help if the condition of Proposition~\ref{isoGM} would also
be necessary for isomorphism. This however is not true!
The isomorphism described in the proposition fixes the switching set $X$
(setwise).
We shall see examples in the next section where $G$ and $G'$ are isomorphic,
but no isomorphism fixes $X$. Because of these examples it will be hard
to find useful conditions for isomorphism that are necessary and sufficient.
Therefore we only present some easy sufficient conditions
for being non-isomorphic after Godsil-McKay switching.
Let $\lambda_G(x,y)$ denote the number of common neighbors of two vertices
$x$ and $y$ in $G$.
It is clear that if the multiset of degrees (i.e.
$\{\lambda_G(x,x)\,|\,x\in V(G)\}$), or the multiset
$\{\lambda_G(x,y)\,|\,x,y\in V(G)\}$ changes after switching,
then $G$ and $G'$ are non-isomorphic. But we can be a bit more precise:

\begin{lem} \label{nonisoGM}
The following conditions are sufficient for $G$ and $G'$ being non-isomorphic.
\begin{itemize}
\item[i] The multiset of degrees (in $G$) of the vertices in $X$ changes
after switching.
\item[ii] The multiset $\Lambda_G=\{\lambda_G(x,y)\,|\,x\in X, y \in V(G)\}$
changes after switching.
\item[iii] The vertices of $X$ all have the same degree, and the multiset
$\overline{\Lambda}_G=\{\lambda_G(x,y)\,|\,x\in X, y \in Y\}$ changes
after switching.
\end{itemize}
\end{lem}
\noindent
{\bf Proof.}
({\em i}) Clearly the degrees in $Y$ don't change by the switching, so the multiset of degrees of $G$ changes whenever the degrees in $X$ change.
({\em ii}) The multiset $\{\lambda_G(x,y)\, |\, x,y \in Y\}$ is not changed after switching, therefore $\{\lambda_G(x,y)\, | x,y\in V(G)\}$ changes if $\Lambda_G(G)$ changes.
({\em iii}) If the vertices in $X$ have the same degree, then switching doesn't change $\{\lambda_G(x,y)\, | x,y\in X\}$.
\qed
\\[5pt]
Suppose not all vertices in $X$ have the same degree.
Then in most cases the set of degrees changes, and hence we get a non-isomorphic graph after switching.
In particular this is always the case if $|X|=4$.

The conditions of Lemma~\ref{nonisoGM} are not necessary for being non-isomorphic.
There are several examples of Godsil-McKay switching in a strongly regular graph $G$
that gives a non-isomorphic graph $G'$ (the smallest example is the
$4\times 4$ grid with a coclique $X$ of size 4).
However, $G'$ is also strongly regular with the same parameters as $G$ (since this property follows from the spectrum), and therefore $\Lambda_G=\Lambda_{G'}$
and $\overline{\Lambda}_G=\overline{\Lambda}_{G'}$.

\section{No isomorphism fixes the switching set}
In this section we give examples of graphs $G$ with a switching set $X$
for which the graphs $G'$ obtained by Godsil-McKay switching are isomorphic
with $G$, but where no isomorphism fixes $X$.

\subsection{Regular tournaments}

A $(0,1)$-matrix $T$ is a {\em tournament matrix} if $T+T^\top=J-I$,
and $T$ is {\em regular} if all row (and column) sums are equal.
If $T$ has order $m$, then this row sum is $(m-1)/2$, so $m$ is odd.

\begin{prp}\label{ex1}
Let $T$ be a regular tournament matrix of order $m > 1$, and put
$N=T\otimes J_2 + I_{2m}$.
Consider a regular graph $H$ of order $2m$ with vertex set $X$ and
automorphism $\rho$ that is a fixed-point-free involution, where
the orbits of the full automorphism group of $H$ are the orbits of $\rho$.
Let $H$ have adjacency matrix $B$, indexed such that $\rho$ is represented
by the permutation matrix $R=I_m\otimes (J_2-I_2)$.
Construct a graph $G$ on the union of two copies $X_1,X_2$ of $X$,
with adjacency matrix
$$
A = \left[\begin{array}{cc}B & N \\
N^{\t} & B \end{array}\right] .
$$
Then $G$ has Godsil-McKay switching set $X_1$, and the switched graph $G'$
is isomorphic with $G$, whilst there is no isomorphism that fixes $X_1$.
\end{prp}
\noindent
{\bf Proof.}
We have $RN=J-N^\top$ and $B = RBR^\top$, and therefore $A'=QAQ^\top$, where
\[
Q = \left[\begin{array}{cc}O & I \\
R & O \end{array}\right].
\]
Thus $G$ is isomorphic with $G'$.
Suppose there is an isomorphism between $G$ and $G'$ that fixes the set
$X_1$ (and hence also $X_2$). Then the isomorphism acts as an automorphism
on the subgraphs induced by $X_1$ and $X_2$, and hence fixes the orbits
of $\rho$ on both copies of $X$. Since $m > 1$ this is impossible.
\qed
\\[5pt]
Regular tournament matrices are easily constructed for every odd order $m$.
If $E$ is the adjacency matrix of an asymmetric regular graph
(asymmetric means that the full automorphism group is trivial),
then $E\otimes J_2$ represents a graph
whose automorphism group satisfies the condition of the proposition.
An asymmetric regular graph exists for every order at least 10
(see~\cite{BI1969}),
but also for $m=5$, 7 and 9 graphs with the required property do exist.
For example when $m=5$ we can take
$$B=\left[\begin{array}{ccccc}
Z & O & Z & O & J\\
O & Z & J & Z & O\\
Z & J & O & Z & O\\
O & Z & Z & O & J\\
J & O & O & J & O
\end{array}\right],
\mbox{ and }
N=\left[\begin{array}{ccccc}
I & J & J & O & O\\
O & I & J & J & O\\
O & O & I & J & J\\
J & O & O & I & J\\
J & J & O & O & I
\end{array}\right],
$$
where $J=J_2$, $I=I_2$ and $Z=J_2-I_2$.
So the construction works for every order $4m$ with $m$ odd and at least 5.
The smallest size of the switching set is 10.
Since in many applications the size of the switching set is 4, the question
rises wether in this special case the sufficient condition for isomorphism
of Proposition~\ref{isoGM} could be necessary.
Unfortunately this is again false, as is illustrated by the next example.

\subsection{A switching set of size four}
Let $G$ be the bipartite
graph on $12+6 = 18$ vertices, where one part of the
bipartition is $\{a,b,c,d,a',b',c',d',a'',b'',c'',d''\}$, the
other is $\{u_i \mid i \in \Zz/6\Zz\}$, and adjacencies are
\begin{eqnarray*}
u_0 \adj a,b,a',c',a'',d'' \\
u_1 \adj b,c,a',b',a'',c'' \\
u_2 \adj b,d,b',c',a'',b'' \\
u_3 \adj c,d,b',d',b'',c'' \\
u_4 \adj a,d,c',d',b'',d'' \\
u_5 \adj a,c,a',d',c'',d'' \\
\end{eqnarray*}
Let the switching set be $X = \{a,b,c,d\}$.
Then we have an isomorphism between $G$ and the switched graph $G'$.
$\phi : G \to G'$ given by $\phi(x) = x'$, $\phi(x') =
x''$, $\phi(x'') = x$ for $x=a,b,c,d$, and $\phi(u_i) =
u_{i+1}$ for $i \in \Zz/6\Zz$.
We would like to show that there is no isomorphism fixing $X$
(but there is). Put $X' = \{a',b',c',d'\}$ and
$X''=\{a'',b'',c'',d''\}$ and $U = \{u_i \mid i \in \Zz/6\Zz\}$.
The graphs $G$ and $G'$ are bipartite and connected,
so any isomorphism $\psi$ fixing $X$ must also fix $X' \cup X''$ and $U$.
The triples $ijk$ such that $u_i,u_j,u_k$ have a common neighbor in $X$ are 045, 012, 135, 234, and after switching 123, 345, 024, 015, so $\psi$ must send the former triples to the latter.
The former triples are precisely the triples with a common neighbor in $X''$, the
latter precisely those with a common neighbor in $X'$.
So $\psi$ must interchange $X'$ and $X''$.
As it turns out, there is such a $\psi$, and we need to enlarge our graph to destroy this unwanted isomorphism.

We can turn the 18-vertex non-example into a 21-vertex almost-example by adding three
vertices $X$, $X'$ and $X''$, corresponding to the sets with the same names, adjacent to their elements (thus: $X \adj a,b,c,d$, etc.), and three directed edges $X \to X'$, $X' \to X''$, and $X'' \to X$.
This gets rid of automorphisms $\psi$ preserving $X$, but the example is directed.
However, Frucht \cite{Frucht1939} showed that every finite
group is the full group of automorphisms of some finite undirected graph.
In particular we can find a graph with full group $C_3$, the cyclic group of order 3, and use that instead of the directed edges.
This yields an actual example.
%
Let us give an explicit example on 9 vertices (\cite{Sab1959}).
Take 9 vertices $x_i$ with $x$ one of $a,b,c$ and $i \in \Zz/3\Zz$.
The 15 edges are $a_i b_i$, $a_i c_{i-1}$, $b_i c_i$, $b_i b_{i+1}$, $b_i c_{i-1}$.
This yields a graph with $C_3$ as full group of automorphisms.
Identify the vertices $X,X',X''$ of the 21-vertex
almost-example with the vertices $a_0$, $a_1$ and $a_2$ of this
gadget (and remove the directed edges) to obtain a 27-vertex
example as claimed.

\section{Graph products}

Consider graphs $G$ and $H$ with adjacency matrices $A$ and $E$, respectively.
We recall that the {\em tensor product} of $H$ and $G$, denoted by $H\times G$
is the graph with adjacency matrix $E\otimes A$.
We will also consider another product, which we will call the {\em strengthened tensor product},
defined by its adjacency matrix $(E+I)\otimes A$, and denoted by $H\stp G$.
Notice that the strengthened tensor product $H\stp G$ can be interpreted as a tensor product
$\overline{H}\times G$ were $\overline{H}$ is obtained from $H$ by adding a loop at every
vertex.

Let $X$ be a switching set in $G$ and suppose that one of the conditions of Lemma~\ref{nonisoGM} is satisfied, so that $G$ is non-isomorphic and cospectral with $G'$.
Then it is easily checked that also the products $G\times H$ and $G\stp H$ are non-isomorphic
and cospectral with $G'\times H$ and $G'\stp H$, respectively.
Indeed, non-isomorphism easily follows because
$\lambda_{H\times G}((i,x),(j,y))=\lambda_{H}(i,j)\lambda_G(x,y)$ and
$\lambda_{H\stp G}((i,x),(j,y))=\lambda_{\overline{H}}(i,j)\lambda_G(x,y)$,
therefore also the multisets
$\{\lambda_{H\times G}((i,x),(j,y))\,|\, i,j\in V(H),\, x,y\in V(G)\}$ and
$\{\lambda_{H\stp   G}((i,x),(j,y))\,|\, i,j\in V(H),\, x,y\in V(G)\}$
are changed after switching (assuming that $H$, resp. $\overline{H}$, has at
least one edge).
Cospectrality follows from basic properties of tensor products of matrices, but also from
the observation that in both products the sets $\{X_i=\{i\}\times X\}$, with $i\in V(H)$,
together with the set $Y$ of remaining vertices is a switching partition.

If none of the conditions of Lemma~\ref{nonisoGM} is satisfied, so that it is conceivable that $G$ is isomorphic with $G'$, then under some easy conditions there exist switching sets in $H\times G$ and $H\stp G$ that lead to non-isomorphic graphs.
For the formulation of the result we will use the notation of Section~2, and the notion of
a pair of {\em complementary rows} in a $(0,1)$-matrix, which simply means that the sum of the two rows is equal to the all-one row.
\begin{thm}\label{prod}
Let $G$ be a graph with a Godsil-McKay switching set $X$, such that the vertices of $X$  have the same degree, and suppose that $\overline{\Lambda}_G=\overline{\Lambda}_{G'}$.
Furthermore suppose that either $X$ is a coclique (i.e. $B=O$), $N$ has at least two columns and no pair of complementary rows,
or that $B$ has row sums $\frac{1}{2}|X|$ and no pair of rows of $\left[\, B\ \, N\,\right]$ is complementary.
Let $H$ be a graph and let $i$ be a vertex of $H$.
Then the subset $\{i\}\times X$ of $V(H) \times V(G)$ is a switching set in $H\times G$
as well as in $H\stp G$, and Godsil-McKay switching gives non-isomorphic cospectral graphs, provided that $i$ has degree at least one in case of the strengthened tensor product and $i$ is adjacent to a vertex of degree at least two in case of the tensor product.
\end{thm}
\noindent
{\bf Proof.}
It is easily checked that for both graph products, the set $\{i\}\times X$
is a switching set.  We'll apply Lemma~\ref{nonisoGM}(iii) and prove that
the multisets $\overline{\Lambda}_{H\times G}$ and
$\overline{\Lambda}_{H\stp G}$ change after switching.

First observe that the Kronecker products $E\otimes A$ and $(E+I)\otimes A$
consist of blocks matrices equal to $A$ or $O$.
After switching the blocks equal to $A$ in the block row and block column
corresponding to $i$ change, but the other blocks remain the same.
For the strengthened tensor product, the diagonal block corresponding to $i$
becomes the switched matrix $A'$.
For both graph products the off-diagonal nonzero blocks in block row $i$ become $A''$, which is obtained
from $A$ by switching with respect to the rows corresponding to $X$.
Note that we can obtain $A''$ also from $A'$ by switching with respect to the
columns corresponding to $X$.
From this it follows that $A''{A''}^\top = A'{A'}^\top$.
\\
For convenience we restrict to the tensor product in the remainder
of the proof; the proof for the strengthened tensor product goes analogously.
The multiset $\overline{\Lambda}_{H\times G}$ consists of the values
$\lambda_{H\times G}((i,x),(j,y))$ where $(i,x) \in \{i\}\times X$
and $(j,y) \not\in \{i\}\times X$.
We distinguish three cases.

\smallskip\noindent
Case (i): $i = j$. We have
\[
\{\lambda_{H\times G}((i,x),(i,y))\,|\,x\in X,y\in Y\}=
\{\lambda_H(i,i)\lambda_G(x,y)\,|\,x\in X,\, y\in Y\},
\]
and $A'{A'}^\top=A''{A''}^\top$ implies that
\[
\{\lambda_{(H\times G)'}((i,x),(i,y))\,|\,x\in X,y\in Y\}=
\{\lambda_H(i,i)\lambda_{G'}(x,y)\,|\,x\in X,\, y\in Y\}.
\]
By assumption the multiset $\overline{\Lambda}_G$ does not change
after switching and therefore the multiset
$\{\lambda_{H\times G}((i,x),(i,y))\,|\,x\in X,y\in Y\}$
is also invariant under switching.

\smallskip\noindent
Case (ii): $i \ne j$ and $y \in Y$. For each $j\neq i$ we have
\[
\{\lambda_{(H\times G)'}((i,x),(j,y))\,|\,x\in X,y\in Y\}=
\{\lambda_H(i,j)\lambda_{G'}(x,y)\,|\,x\in X,\, y\in Y\}=
\]
\[
\{\lambda_H(i,j)\lambda_G(x,y)\,|\,x\in X,\, y\in Y\}=
\{\lambda_{H\times G}((i,x),(j,y))\,|\,x\in X,y\in Y\}.
\]
Case (iii): $i \ne j$ and $x,y \in X$.
Choose $\j\neq i$ such that $\lambda_H(i,\j)$ is maximal.
It follows that $\lambda_H(i,\j)>0$ because $i$ has a neighbor of degree
at least two. (Note that for the strengthened tensor product it suffices
that the degree of $i$ is at least 1.)
We have $\lambda_{H\times G}((i,x),(\j,x))=\lambda_H(i,\j)\lambda_G(x,x)$.
After switching we get
$\lambda_{(H\times G)'}((i,x),(\j,x))=\lambda_H(i,\j)\mu(x)$,
where $\mu(x)$ is the number of neighbors of $x$ that remain a neighbor
after switching. Clearly $\mu(x)<\lambda_G(x,x)$, hence
$$
\lambda_{(H\times G)'}((i,x),(\j,x)) <
\lambda_{H\times G}((i,x),(\j,x)) .
$$
For $y \ne x$ we get
$\lambda_{(H\times G)'}((i,x),(j,y))=\lambda_H(i,j)\lambda_{G'}(x,y)$.
Because the matrices $N$ or $[\, B\ N\,]$ which are switched to their
complements have no complementary pair of rows, it follows that
$\lambda_{G'}(x,y)<\lambda_G(x,x)$.
Hence we have
\[
\lambda_{(H\times G)'}((i,x),(j,y)) < \lambda_H(i,\j)\lambda_{G}(x,x)
= \lambda_{H\times G}((i,x),(\j,x)).
\]
This implies that the number $\lambda_{H\times G}((i,x)(\j,x))$
disappears at least once from the multiset $\overline{\Lambda}_{H\times G}$
after switching.
\qed
\\[5pt]
In view of the previous section it seems relevant to remark that the proof
of the above theorem would have been much simpler if we could have used
that there exists an isomorphism that fixes the switching set.

The $\ell\times m$ grid (or lattice graph $L(\ell,m)$) is the line graph of the complete
bipartite graph $K_{\ell,m}$ (we assume $\ell\geq m$).
If $(\ell,m)\neq(4,4)$, or $(6,3)$, then $L(\ell,m)$ is determined by its spectrum.
If $\ell>m\geq 2$ a 4-cycle in the grid is a switching set that satisfies the hypothesis of Theorem~\ref{prod}.
Therefore the tensor product of $L(\ell,m)$ ($\ell>m\geq 2$) and a graph with
at least one vertex of degree two is not determined by its adjacency spectrum.

The strengthened tensor product $K_n\stp G$ ($n>1$) is also known
as a {\em coclique extension} of $G$.
So the above theorem gives some easy conditions for a coclique extension
to have non-isomorphic cospectral graphs.
For example a coclique extension of the grid $L(\ell,m)$
with $\ell>m\geq 2$, is not determined by its spectrum.

Another example is the triangular graph $T(m)$, which is the line graph of $K_m$.
If $m\neq 8$ the spectrum determines $T(m)$ and if $m\geq 4$ a 4-cycle in $T(m)$
satisfies the requirements of Theorem~\ref{prod}.
Thus we can conclude that for $m\geq 4$ a coclique extension of $T(m)$ is not determined by its spectrum.

\end{document}